\documentclass[12pt]{amsart}
\usepackage{amscd}     
\usepackage{amssymb}
\usepackage{amsmath, amsthm, graphics}
\usepackage{xypic}     
\LaTeXdiagrams        
\usepackage[all]{xy}
\xyoption{2cell} \UseAllTwocells \xyoption{frame} \CompileMatrices
\allowdisplaybreaks[3]
\usepackage{amsfonts}
\usepackage[normalem]{ulem}
\usepackage{hyperref}
\usepackage{tikz}
\usepackage{mathtools}
\usepackage{verbatim}

\usepackage{latexsym}
\usepackage{epsfig}
\usepackage{amsfonts}
\usepackage{enumerate}
\usepackage{times}

\newtheorem{prop}{Proposition}

\newtheorem{question}[prop]{Question}

\newtheorem{theorem}{Theorem}[section]

\theoremstyle{remark}

\theoremstyle{remark}

\newtheorem{remark}[subsection]{Remark}

\numberwithin{equation}{section}

\newcommand{\Mbar}{\overline{\M}}

\newcommand{\X}{\mathcal{X}}
\newcommand{\Y}{\mathcal{Y}}

\newcommand{\K}{\mathcal{K}}
\newcommand{\M}{\mathcal{M}}
\newcommand{\C}{\mathcal{C}}

\newcommand{\D}{\mathcal{D}}

\def\<{\left\langle}
\def\>{\right\rangle}

\def\b1{{\mathbf 1}}

\begin{document}

\title{On orbifold Gromov-Witten classes}

\author[Tseng]{Hsian-Hua Tseng}
\address{Department of Mathematics\\ Ohio State University\\ 100 Math Tower, 231 West 18th Ave. \\ Columbus,  OH 43210\\ USA}
\email{hhtseng@math.ohio-state.edu}

\date{\today}

\maketitle

{\centering \em Dedicated to the loving memory of Ching-Shan Chou\par}


\setcounter{tocdepth}{1}

\setcounter{section}{-1}
\section{Introduction}
Let $\X$ be a smooth proper Deligne-Mumford stack over $\mathbb{C}$. The stack 
\begin{equation*}
    \mathcal{K}_{g,n}(\X,d),
\end{equation*}
which parametrizes degree $d$ stable maps from genus $g$ orbifold curves\footnote{Orbifold curves are also called {\em twisted curves}. Orbifold nodes are always assumed to be {\em balanced}.} with $n$ possibly orbifold markings\footnote{The marked gerbes are {\em not} trivialized.}, is constructed in \cite{av} (see also \cite{ol}). It is Deligne-Mumford and proper over $\mathbb{C}$. 

There are several natural maps defined for $\mathcal{K}_{g,n}(\X,d)$:
\begin{enumerate}
    \item 
    Restricting stable maps to marked points yields the {\em evaluation maps} 
    \begin{equation*}
        ev: \mathcal{K}_{g,n}(\X,d)\to \bar{I}\X,
    \end{equation*}
    where $\bar{I}\X$ is the {\em rigidified} inertia stack of $\X$. See \cite[Section 3]{agv} for a detailed discussion on inertia stacks and \cite[Section 4.4]{agv} for the construction of evaluation maps.
    \item
    Forgetting stable maps to $\X$ but only retaining the domain curves yields the forgetful map 
    \begin{equation*}
        \pi:\mathcal{K}_{g,n}(\X,d)\to \mathfrak{M}_{g,n}^{tw},
    \end{equation*}
    where $\mathfrak{M}_{g,n}^{tw}$ is the stack of $n$-pointed genus $g$ orbifold curves, see \cite[Theorem 1.9]{ol}. Assuming $2g-2+n>0$, then passing to coarse curves and stablizing the domains yield another forgetful map
    \begin{equation*}
        p: \mathcal{K}_{g,n}(\X,d)\to \overline{\M}_{g,n},
    \end{equation*}
    where $\overline{\M}_{g,n}$ is the stack of $n$-pointed genus $g$ stable curves. There is an obvious commutative diagram
    \begin{displaymath}
    \xymatrix{ 
    \mathcal{K}_{g,n}(\X,d)\ar[rd]^{p}\ar[d]_{\pi} &  \\
   \mathfrak{M}_{g,n}^{tw}\ar[r]^{}& \overline{\M}_{g,n}.}
\end{displaymath}
\end{enumerate}

A perfect obstruction theory for $\mathcal{K}_{g,n}(\X, d)$ relative to $\pi$ is introduced in \cite{agv}, yielding a virtual fundamental class in Chow groups\footnote{Chow and (co)homology groups are taken with $\mathbb{Q}$-coefficients.},
\begin{equation*}
    [\mathcal{K}_{g,n}(\X,d)]^{vir}\in \text{CH}_*(\mathcal{K}_{g,n}(\X,d)),
\end{equation*}
 which may also be viewed as a homology class via the cycle map $\text{CH}_*(\mathcal{K}_{g,n}(\X,d))\to H_*(\mathcal{K}_{g,n}(\X,d))$.

There are natural classes defined on $\mathcal{K}_{g,n}(\X,d)$:
\begin{enumerate}
    \item 
    Pulling back via evaluation maps yields $$ev_i^*(\gamma)$$ where $\gamma$ is a Chow/cohomology class of $\bar{I}\X$. 
    \item
    The {\em descendant classes} $$\psi_i:=c_1(L_i)$$ are the first Chern classes (taken in Chow or cohomology groups) of line bundles $L_i\to \mathcal{K}_{g,n}(\X,d)$ formed by cotangent lines at the $i$-th marked points of the {\em coarse domain curves}.
\end{enumerate}
Goromov-Witten theory of the stack $\X$ is the study of classes of the form
\begin{equation}\label{eqn:oGW_class1}
    \prod_{i=1}^n\psi_i^{k_i} ev_i^*(\gamma_i)\cap [\mathcal{K}_{g,n}(\X,d)]^{vir},
\end{equation}
where $\gamma_1,...,\gamma_n$ are Chow/cohomology classes of $\bar{I}\X$ and $k_1,...,k_n\in \mathbb{Z}_{\geq 0}$. Pushing (\ref{eqn:oGW_class1}) to a point yields Gromov-Witten {\em invariants} of $\X$, which have been studied extensively in the past $20$ years. The purpose of this note is to discuss some questions arising from pushing forward (\ref{eqn:oGW_class1}) to other natural settings.

\subsection*{Acknowledgement}
We thank R. Pandharipande, J. Schmitt, and the referee for helpful comments. The author is supported in part by Simons Foundation Collaboration Grant.

\section{Tautological cohomology classes}\label{sec:taut_class}
In (\ref{eqn:oGW_class1}), take $\gamma_1,...,\gamma_n\in H^*(\bar{I}\X)$ to be cohomology classes. Pushing forward (\ref{eqn:oGW_class1}) via $p$ yields what is usually called {\em Gromov-Witten classes}:
\begin{equation}\label{eqn:oGW_class2}
    p_*\left(\prod_{i=1}^n\psi_i^{k_i} ev_i^*(\gamma_i)\cap [\mathcal{K}_{g,n}(\X,d)]^{vir}\right)\in H^*(\overline{\M}_{g,n}).
\end{equation}
Without descendants, the classes (\ref{eqn:oGW_class2}) yield a system of multi-linear maps 
\begin{equation}\label{eqn:GW_CohFT}
    H^*(\bar{I}\X)^{\otimes n}\to H^*(\overline{\M}_{g,n}), \quad \gamma_1\otimes...\otimes\gamma_n\mapsto p_*\left(\prod_{i=1}^nev_i^*(\gamma_i)\cap [\mathcal{K}_{g,n}(\X,d)]^{vir}\right).
\end{equation}
Properties of virtual fundamental classes imply that (\ref{eqn:GW_CohFT}) is a {\em cohomological field theory}, a notion introduced in \cite{km}. Discussions on further developments of cohomological field theories can be found in \cite{p}.

An important aspect of the study of $H^*(\overline{\M}_{g,n})$ is the (cohomological) {\em tautological ring}
\begin{equation*}
    RH^*(\overline{\M}_{g,n})\subset H^*(\overline{\M}_{g,n}),
\end{equation*}
which can be defined as the smallest system of unital subrings of $H^*(\overline{\M}_{g,n})$ which is stable under push-forward and pull-back by the following maps:
\begin{enumerate}
    \item 
    $\Mbar_{g,n+1}\to \Mbar_{g,n}$ forgetting one of the markings;
    \item
    $\Mbar_{g_1,n+1}\times \Mbar_{g_2,n_2+1}\to \Mbar_{g_1+g_2, n_1+n_2}$ gluing two curves at a point;
    \item
    $\Mbar_{g-1, n+2}\to \Mbar_{g,n}$ gluing together two points on a curve.
\end{enumerate}
More details can be found in e.g. \cite{fp}.

Elements of $RH^*(\Mbar_{g,n})$ are called {\em tautological classes}. The following question, raised for smooth projective varieties \cite{fp}, should obviously be asked for stacks:

\begin{question}\label{Q:tautological}
Let $\X$ be a smooth proper Deligne-Mumford stack over $\mathbb{C}$. For $\gamma_1,...,\gamma_n\in H^*(\bar{I}\X)$ and $k_1,...,k_n\in \mathbb{Z}_{\geq 0}$, are the Gromov-Witten classes (\ref{eqn:oGW_class2}) tautological?
\end{question}
\begin{remark}
While tautological rings inside the Chow ring $\text{CH}^*(\Mbar_{g,n})$ can be defined, the Chow version of Question \ref{Q:tautological} is not expected to be true even for varieties. Hence we do not discuss the Chow version here.
\end{remark}

Question \ref{Q:tautological} is known to be true for a number of classes of varieties. A summary of known results can be found in \cite[Section 0.7]{abpz}. Here we provide two classes of Deligne-Mumford stacks for which Question \ref{Q:tautological} is true.

\begin{theorem}\label{thm:toric_taut}
Question \ref{Q:tautological} is true for smooth semi-projective {\em toric} Deligne-Mumford stacks $\X$.
\end{theorem}
Just like the case of toric varieties, Theorem \ref{thm:toric_taut} follows from virtual localization \cite{gp}. Virtual localization formula for Gromov-Witten theory of toric Deligne-Mumford stacks is written very explicitly in \cite{l}. 

Virtual localization formula reduces Theorem \ref{thm:toric_taut} to studying the {\em Hurwitz-Hodge classes}. To address this, we first review the construction of Hurwitz-Hodge classes arising in the present setting. Let $G$ be a finite abelian group. Let $V$ be a finite dimensional $\mathbb{C}$-vector space that admits a $G$-action. Let $T$ be an algebraic torus with an action on $V$ (so that the $G$ and $T$ actions on $V$ commute). The vector space $V$ defines a $T$-equivariant vector bundle $\mathcal{V}\to BG$. Let $\mathcal{K}_{g,n}(BG)$ be the moduli stack of stable maps to $BG=[\text{pt}/G]$. Consider the universal stable map,
\begin{displaymath}
    \xymatrix{ 
    \mathcal{C}\ar[r]^{f}\ar[d]_{q} & BG \\
   \mathcal{K}_{g,n}(BG).& }
\end{displaymath}
The $K$-theory class $$Rq_*f^*\mathcal{V}\in K_*(\mathcal{K}_{g,n}(BG))$$
is well-defined. Hurwitz-Hodge classes are $T$-equivariant inverse Euler classes $e_T^{-1}(Rq_*f^*\mathcal{V})$ of this kind of $K$-theory objects.

Now, $e_T^{-1}(Rq_*f^*\mathcal{V})$ can be expressed in terms of Chern characters of $Rq_*f^*\mathcal{V}$. The (more general) Riemann-Roch calculation\footnote{Strictly speaking, the Riemann-Roch calculation in \cite{t} is done for a different moduli stack $\Mbar_{g,n}(BG)$ parametrizing stable maps {\em with sections to marked gerbes}. The answer can be easily adjusted to the present setting.} of \cite{t} implies that these Chern characters can be expressed in terms of $\psi$ classes and boundary classes of $\mathcal{K}_{g,n}(BG)$. Hence, after pushing forward to $\Mbar_{g,n}$, Hurwitz-Hodge classes $e_T^{-1}(Rq_*f^*\mathcal{V})$ are tautological. This proves Theorem \ref{thm:toric_taut}.

\begin{remark}
\hfill
\begin{enumerate}
    \item 
    The above argument is valid for Deligne-Mumford stacks $\X$ admitting torus actions with isolated fixed points and $1$-dimensional orbits.
    \item 
    The above argument is valid in Chow groups, thus answering the Chow version of Question \ref{Q:tautological} in affirmative for toric $\X$.
    \item
    Virtual localization was applied to study Gromov-Witten theory of {\em toric bundles} $E\to B$ in \cite{cgt}. It is clear from the localization analysis in \cite{cgt} and Riemann-Roch calculations in \cite{fp_hh} and \cite{cg} that Gromov-Witten classes of the toric bundle $E$ are tautological if Gromov-Witten classes of the base $B$ are tautological. This gives another evidence for Question \ref{Q:tautological}.
\end{enumerate}
\end{remark}

The second class of examples we consider is orbifold curves. It is known that a smooth orbifold curve $\C$ is obtained from its underlying coarse curve $C$ (which is itself a smooth curve) by applying a finite number of {\em root constructions}. We refer to \cite[Theorem 4.2.1]{agv} for more details of this description. 

\begin{theorem}\label{thm:curve_taut}
Question \ref{Q:tautological} is true for smooth projective orbifold curves $\C$.
\end{theorem}

Question \ref{Q:tautological} is proven to be true for nonsingular curves in \cite{j}. Our proof of Theorem \ref{thm:curve_taut} builds on that result, as follows.

We begin with some notations. Let $p_1,...,p_m\in \C$ be the orbifold points of $\C$, and $\bar{p}_1,...,\bar{p}_m\in C$ their images in the coarse curve. Note that $\bar{p}_1,...,\bar{p}_m$ are smooth points on $C$. Let $r_1,...,r_m\in \mathbb{N}$ be orders of stabilizers of the orbifold points $p_1,...,p_m$ respectively. Deformation to the normal cone construction can be applied to $p_1,...,p_m\in \C$ to give a degeneration of $\C$ to the following nodal curve:
\begin{equation}\label{eqn:nodal_curve}
    C\bigcup\cup_{i=1}^m \mathbb{P}_{1, r_i}^1,
\end{equation}
where $\bar{p}_i\in C$ is identified with the smooth point $0\in \mathbb{P}_{1,r_i}^1$. More precisely, this degeneration is obtained by degenerating the coarse curve $C$ to $C\bigcup\cup_{i=1}^m \mathbb{P}^1$ (where $\bar{p}_i$ is identified with $0\in \mathbb{P}^1$), then applying the $r_i$-th root construction to the divisor in the total space formed as $\bar{p}_i$ moves.

Associated to the pairs $(C,\bar{p}_1,...,\bar{p}_m)$ and $\{(\mathbb{P}_{1,r_i}^1, 0)\}_{i=1}^m$ are their {\em relative} Gromov-Witten classes. Relative Gromov-Witten classes of a pair $(\X, \D)$ of a smooth proper  Deligne-Mumford stack $\X$ and a smooth divisor $\D\subset \X$ are defined in a manner similar to (\ref{eqn:oGW_class2}) by working with moduli stacks of stable relative maps to $(\X, \D)$. Details of these moduli stacks can be found in \cite{af}.

Degeneration formula, proven in \cite{af}, applies to this setting and expresses Gromov-Witten classes (\ref{eqn:oGW_class2}) of $\C$ in terms of relative Gromov-Witten classes of $(C,\bar{p}_1,...,\bar{p}_m)$ and $\{(\mathbb{P}_{1,r_i}^1, 0)\}_{i=1}^m$. By \cite[Theorem 1]{j}, relative Gromov-Witten classes of $(C,\bar{p}_1,...,\bar{p}_m)$ are tautological. The pair $(\mathbb{P}_{1,r_i}^1,0)$ is toric, and the relative virtual localization formula may be applied. The argument described in the proof of Theorem \ref{thm:toric_taut} applies here to show that some terms in relative virtual localization formula are tautological. The only terms not covered by this argument are the {\em double ramification cycles}, which are tautological by \cite{fp} or \cite{jppz1}. Therefore relative Gromov-Witten classes of $(\mathbb{P}_{1,r_i}^1,0)$ are tautological. This proves Theorem \ref{thm:curve_taut}.

\begin{remark}
The above argument can be extended a little bit to show the following: for a smooth projective variety $X$ and a smooth divisor $D\subset X$, Gromov-Witten classes of the stack $X_{D,r}$ of $r$-th roots of $X$ along $D$ are tautological if relative Gromov-Witten classes of $(X,D)$ and absolute Gromov-Witten classes of $D$ are tautological. This proof encounters double ramification cycles with target $D$, which are tautological by the formula in \cite{jppz2}, provided that Gromov-Witten classes of $D$ are tautological.
\end{remark}

\begin{remark}
It would be interesting to consider Question \ref{Q:tautological} in other examples. For instance, it follows from the product formula \cite{ajt} that Gromov-Witten classes of a product stack $\X\times \Y$ are tautological if Gromov-Witten classes of $\X$ and $\Y$ are tautological. With efforts, one can hope that the approach in \cite{abpz} can be extended to complete intersections in weighted projective stacks.
\end{remark}

\section{Global finite group quotients}
An important aspect of the Gromov-Witten theory of stacks $\X$ is the presence of orbifold structures in the domains of stable maps to $\X$. The morphism $p:\K_{g,n}(\X,d)\to \Mbar_{g,n}$ forgets these orbifold structures. Therefore it is interesting to consider Gromov-Witten classes of $\X$ in suitable settings where these orbifold structures are not forgotten. 

Here we discuss an attempt to retain these orbifold structures for target stacks of the form
\begin{equation*}
    \X=[M/G],
\end{equation*}
where $M$ is a smooth (quasi)projective variety over $\mathbb{C}$ and $G$ is a {\em finite} group. The constant map $M\to \text{pt}$ is clearly $G$-equivariant, and yields a representable morphism 
\begin{equation*}
    \X=[M/G]\to BG=[\text{pt}/G].
\end{equation*}
Composing stable maps to $\X$ with this morphism and stablizing yield a morphism of moduli stacks
\begin{equation*}
    p_G:\K_{g,n}(\X,d)\to\K_{g,n}(BG),
\end{equation*}
which is proper. Pushing forward (\ref{eqn:oGW_class1}) via $p_G$ yields the following classes   
\begin{equation}\label{eqn:oGW_class3}
    (p_G)_*\left(\prod_{i=1}^n\psi_i^{k_i} ev_i^*(\gamma_i)\cap [\mathcal{K}_{g,n}(\X,d)]^{vir}\right)\in H^*(\K_{g,n}(BG)).
\end{equation}
Further pushing forward (\ref{eqn:oGW_class3}) via the natural map $\K_{g,n}(BG)\to \overline{\M}_{g,n}$ recovers (\ref{eqn:oGW_class2}).

An interesting subring of $H^*(\K_{g,n}(BG))$ is the (cohomological) {\em $\mathcal{H}$-tautological ring}\footnote{It is originally defined in the Chow theory.}:
\begin{equation*}
    R_{\mathcal{H}}(\K_{g,n}(BG))\subset H^*(\K_{g,n}(BG)),
\end{equation*}
see \cite{lian}.
\begin{question}\label{Q:tautological2}
Are (\ref{eqn:oGW_class3}) contained in the $\mathcal{H}$-tautological ring of $\K_{g,n}(BG)$?
\end{question}
When $M$ is toric, and the $G$-action commutes with the torus action on $M$, the stack $\X=[M/G]$ admits a torus action and the above approach to Theorem \ref{thm:toric_taut} applies to show that Question \ref{Q:tautological2} is true in this case. However, in this case (\ref{eqn:oGW_class3}) are contained in some smaller subset of $R_{\mathcal{H}}(\K_{g,n}(BG))$. Indeed, virtual localization formula and Riemann-Roch calculations show that (\ref{eqn:oGW_class3}) are obtained from pushforwards of combinations of $\psi$ classes via the following natural morphisms: 
\begin{enumerate}
    \item The morphism that forgets a non-stacky marking, as discussed in \cite[Proposition 8.1.1]{agv};
    
    \item The boundary gluing morphisms, as discussed in \cite[Proposition 5.2.1]{agv}.
\end{enumerate}
It should be possible to define a subring of $H^*(\K_{g,n}(BG))$ using the definition of $RH^*(\Mbar_{g,n})$, recalled in Section \ref{sec:taut_class}, with these maps. If defined, this subring is smaller than the $\mathcal{H}$-tautological ring. Still, (\ref{eqn:oGW_class3}) lie in such a subring.

Whether the formulation of Question \ref{Q:tautological2} really requires the $\mathcal{H}$-tautological rings remains unclear.

For more general $\X$, it is not clear how to construct variants of Gromov-Witten classes of $\X$ that retain orbifold structures on the domains. The natural place for keeping the domain orbifold curves is the stack $\mathfrak{M}_{g,n}^{tw}$ of orbifold curves. However the morphism $\pi: \K_{g,n}(\X,d)\to \mathfrak{M}_{g,n}^{tw}$ is not necessarily proper and cannot be used to produce interesting classes on $\mathfrak{M}_{g,n}^{tw}$, even though the tautological Chow ring $R^*(\mathfrak{M}_{g,n}^{tw})$ can be defined\footnote{With known results in \cite{acv} and \cite{ol}, it is not hard to see that the approach in \cite{bs} can be adopted to define $R^*(\mathfrak{M}_{g,n}^{tw})$, see \cite{t22}.}.

\end{document}